\documentclass[12pt,leqno]{article}
\usepackage[pdftex]{graphicx}
\usepackage{amsfonts, amsmath}
\usepackage{times,wrapfig,verbatim}
\usepackage{color}

\newcounter{conjecture}
\newcounter{remark}
\newcounter{corollary}

\newenvironment{corollary}{\medskip{\bf Corollary \thecorollary.}
\addtocounter{corollary}{1}\em}{\rm}
\newtheorem{theorem}{Theorem}
\newtheorem{lemma}{Lemma}
\newtheorem{proposition}{Proposition}
\newtheorem{conject}{Conjecture}

\newcommand{\dd}{\delta}

\newcommand{\lar}{\longrightarrow}
\newcommand{\eps}{\varepsilon}

\newcommand{\lll}{\label}
\newcommand {\rrr}[1]{(\ref{#1})}

\def \be{\begin{equation}}
\def \ee{\end{equation}}
\def \bt{\begin{theorem}}
\def \et{\end{theorem}}
\def \bc{\begin{corollary}}
\def \ec{\end{corollary}}
\def \bea{\begin{eqnarray}}
\def \eea{\end{eqnarray}}
\def \bas{\begin{eqnarray*}}
\def \eas{\end{eqnarray*}}

\def \aa{\alpha}

\def \dd{\delta}

\def \si{\sigma}

\def \th{\theta}

\def \vski{\vspace{12pt}}

\def \ff{\infty}

\def \AA{{\cal A}}

\def \({\left(}
\def \){\right)}

\def \bc{\begin{center} }
\def \ec{\end{center} }
\def \bs{\begin{slide} }
\def \es{\end{slide} }

\def\square{{\vcenter{\vbox{\hrule height.3pt
         \hbox{\vrule width.3pt height5pt \kern5pt
            \vrule width.3pt}
         \hrule height.3pt}}}}
\def\qed{{\hfill $\square$ \bigskip}}

\newcounter{cccases}
\numberwithin{equation}{section}

\newcommand {\SC} {{\mathbb C}}
\newcommand {\SD} {{\mathbb D}}

\newcommand {\SR} {{\mathbb R}}

\begin{document}

\title{On the probability of fast exits and long stays of planar Brownian motion in simply connected domains}

\author{Dimitrios Betsakos$^{\,\rm 1}$,~Maher Boudabra$^{\,\rm 2}$, and Greg Markowsky$^{\,\rm 2}$\\
{\small {\tt betsakos@math.auth.gr} ~~ {\tt maher.boudabra@monash.edu} ~~ {\tt gmarkowsky@gmail.com}}\\
{\footnotesize{$^{\rm 1}$Department of Mathematics,  Aristotle University of Thessaloniki, Greece}}\\
{\footnotesize{$^{\rm 2}$Department of Mathematics,  Monash University, Australia}}}




\maketitle

\begin{abstract}
Let $T^D$ denote the first exit time of a planar Brownian motion
from a domain $D$. Given two simply connected planar domains $U,W \neq \SC$ containing $0$,
we investigate the cases in which we are more likely to have fast exits
(meaning for instance ${\bf P}(T^U<t) > {\bf P}(T^W<t)$ for $t$ small) from $U$ than from $W$,
or long stays (meaning ${\bf P}(T^U>t) > {\bf P}(T^W>t)$ for $t$ large).
We prove several results on these questions. In particular, we show that the primary factor in the probability of fast exits is the proximity of the
boundary to the origin, while for long stays an important factor is the moments of the exit time. The complex analytic theory that motivated our inquiry is also discussed.
\end{abstract}


\section{Introduction}
\label{1}

The distribution of  the exit time of planar Brownian motion from
a domain measures in some sense the size of the domain. It can
also be used for the study of analytic functions via the conformal
invariance of Brownian motion.

\medskip

Let $Z_t=X_t+iY_t$, $t\geq 0$ denote standard Brownian motion
moving in the plane and starting from the origin, (that is,
$Z_0=0$ almost surely). We denote by ${\bf P}$ and ${\bf E}$ the
corresponding probability measure and expectation, respectively.
For a domain $D$ containing $0$ in the complex plane $\SC$, we denote by $T^D$
the first exit time of $Z_t$ from $D$; that is
$$
T^D=\inf\{t>0: Z_t\notin D\}.
$$

\medskip

Suppose also that $f$ is a univalent function in the well-known
class $\mathcal S$. Thus $f$ is a function univalent and
holomorphic in the unit disk $\SD$ with $f(0)=0$ and
$f^\prime(0)=1$. By a classical result of P. L\'{e}vy, the image
of $Z_t$ under $f$ is a Brownian motion with a time change. We
describe now a precise version of this conformal invariance of
Brownian motion; (see \cite[\S 2]{DAVIS}). Let
$$
\rho_f(s)=\int_0^s |f^\prime(Z_t)|^2\,dt,\;\;\;\;\; 0\leq s<T^\SD.
$$
Observe that $\rho_f$ is almost surely strictly increasing and set
$$
W_t=f(Z_{\rho_f^{-1}(t)}),\;\;\;0\leq t<\rho_f(T^\SD).
$$
Define also $W_{\rho_f(T^\SD)}=\lim_{t\to T^\SD}W_{\rho_f(t)}$ and
$$
W_{\rho_f(T^\SD)+t}=W_{\rho_f(T^\SD)}+(Z_{T^\SD+t}-Z_{T^\SD}),\;\;\;\;t>0.
$$
Note that $f$ being univalent and holomorphic implies that $P(\rho_f(T^\SD)< \ff) = 1$; this is a nontrivial statement, but is shown for instance in \cite[p.198]{burk}. L\'{e}vy's theorem now asserts that $\{W_t,t\geq 0\}$ is standard
planar Brownian motion starting from the origin.

\medskip

We set
$$
\nu(f)=\rho_f(T^\SD)=\int_0^{T^\SD}|f^\prime(Z_t)|^2\,dt
$$
and observe that $\nu(f)$ is the first exit time of $W_t$ from
$f(\SD)$; that is, $\nu(f)=T^{f(\SD)}$. ``The distribution of
$\nu(f)$ is an intuitively appealing measure of the size of
$f(\SD)$" \cite{DAVIS}. In several classical extremal problems for
functions in the class ${\mathcal S}$, the identity function
$I(z)=z$ is the ``smallest" function while the Koebe function
$k(z)=z/(1-z)^2$ is the ``largest" function in $\mathcal S$. B.
Davis \cite{DAVIS} conjectured that \be\label{i1} {\bf
E}[\Phi(\nu(I))]\leq {\bf E}[\Phi(\nu(f))], \ee for all
$f\in\mathcal S$ and all increasing convex functions
$\Phi:[0,\infty)\to\SR$, and suggested that perhaps
$$
{\bf P}(\nu(I)>t) \leq {\bf P}(\nu(g)>t),
$$
for all $t>0$ and  $g\in\mathcal S$. Apropos the  first conjecture
T. McConnell \cite{MAC} proved that \be\label{i1.5} {\bf
E}[\nu(I)^p]\leq {\bf E}[\nu(f)^p],\;\;\;0<p<\infty, \ee but the
full conjecture remains open, as far as we know. McConnell also
disproved the second conjecture by finding functions $g\in
\mathcal S$ such that \be\label{i2} {\bf P}(\nu(I)>t)>{\bf
P}(\nu(g)>t), \ee for all sufficiently small $t>0$. Davis also
asked in what sense, with regard to $\nu(f)$, the Koebe function
is the largest in $\mathcal S$.

\medskip


Motivated by these developments, we have considered the following questions.
Given two simply connected planar domains $U,W \neq \SC$ containing $0$,
what sufficient conditions can we place on the domains so that we are more likely to have fast exits
(meaning for instance ${\bf P}(T^U<t) > {\bf P}(T^W<t)$ for $t$ small) from $U$ than from $W$,
or long stays (meaning ${\bf P}(T^U>t) > {\bf P}(T^W>t)$ for $t$ large).
We have found that the primary factor influencing the probability of fast exits is the proximity of the
boundary to the origin. In order to make this more precise, let us introduce the following notation. For any simply connected domain $V$,
let
$$
d(V) = \inf\{|z|: z \in \partial V\}.
$$

We then have the following theorem, which is the main result of the paper.

\begin{theorem}\label{T2}
Suppose that $d(U) < \frac{1}{\sqrt{2}} d(W)$. Then, for all sufficiently small $t>0$,

\begin{equation}\label{T2e}
{\bf P}(T^U<t)> {\bf P}(T^W< t).
\end{equation}

In fact,

$$
\lim_{t \lar 0^+} \frac{{\bf P}(T^U<t)}{{\bf P}(T^W< t)} = \ff.
$$
\end{theorem}

We do not know whether $\frac{1}{\sqrt{2}}$ is the optimal constant; it may even be that it may be replaced by 1. If so, this would be a surprising property of planar Brownian motion. Evidence for this possibility can be found in the fact that it is true if $U$ is a half-plane. This is discussed in more detail in Section \ref{conc}.

\vski

Naturally, it would be nice to have an analog for long stays. We believe that the important factor for long stays in domains is the moments of the exit time. To be precise, for domain $V$ let
$$
{\rm H}(V)=\sup \{p > 0: {\bf E}[(T^V)^p] < \ff\};
$$
note that ${\rm H}(V)$ is proved in \cite{burk} to be exactly
equal to half of the Hardy number of $V$, a purely analytic
quantity, as defined in \cite{hansen}, and is therefore calculable for a number of common domains. Furthermore ${\rm H}(V)
\geq \frac{1}{4}$ as long as $V \neq \SC$. We have the following simple result.

\begin{proposition}\label{T3}
Suppose that ${\rm H}(U) > {\rm H}(W)$. Then
\begin{equation}\label{T3e}
\limsup_{t \lar \ff} \frac{{\bf P}(T^W > t)}{{\bf P}(T^U > t)} = \ff.
\end{equation}
\end{proposition}

We conjecture that this proposition is true with the $\limsup$ replaced by $\lim$, but have not been able to prove it (except when $W$ is either a half-plane or quarter-plane). This is discussed in detail in the final section.




\section{Proofs}
\label{2}

In this section we prove Theorem \ref{T2} and Proposition \ref{T3}.

\subsection{Proof of Theorem \ref{T2}}

The proof of Theorem \ref{T2} is based on the strong Markov
property and the explicit formula for the transition density of
one-dimensional Brownian motion. Some of the estimates used are
probably known to experts; we hope, however, that their elementary
derivation and their use in the study of univalent functions are
of some interest.

\vski

In what follows, $X_t$ will denote one-dimensional Brownian motion.
The corresponding probability measure with starting point $x\in
\SR$ will be denoted by ${\bf P}^x$. The first hitting time of a
point $y\in\SR$ will be denoted by $\tau_y$. The following well
known equality comes easily from the reflection principle; see

\cite[p.23]{DurBM}:
\begin{equation}\label{p1}
{\bf P}^0(\tau_a\leq t)=2{\bf P}^0(X_t\geq a),\;\;\;\;a>0.
\end{equation}

We will use the standard notation for the transition density
function of Brownian motion:
$$
p_s(x)=\frac{1}{\sqrt{2\pi s}}\;\;e^{-x^2/2s},\;\;\; x\in\SR,
\;s>0.
$$
It follows from elementary calculus that, for any $\dd > 0$ there
exists a constant $C_1>1$ such that for every $y>\dd$,

\begin{equation}\label{L1e}
C_1^{-1} \;\frac{e^{-y^2}}{y}\leq \int_y^\infty e^{-x^2}\,dx\leq
C_1 \;\frac{e^{-y^2}}{y}.
\end{equation}

We begin by proving a preliminary proposition, then show how it extends to prove Theorem \ref{T2}.

\begin{proposition} \label{prelim}
Let $K_\aa := \SC \backslash (-\ff,-\aa]$ for any $\aa$ with $0<\aa<\frac{1}{\sqrt{2}}$. Then

$$
\lim_{t \lar 0^+} \frac{{\bf P}(T^{K_\aa}<t)}{{\bf P}(T^\SD< t)} = \ff.
$$

\end{proposition}

This will be proved through a sequence of lemmas. In what follows, $C$ will denote a generic absolute constant that
may change from line to line.

\begin{lemma}\label{L2}
For any $\eps>0$, there exist constants $\dd,C_2>0$ such that for every $t\in
(0,\dd)$,
\begin{equation}\label{L2e}
{\bf P}^0(X_{t/2}>0,X_t<0) = \int_0^\infty p_{t/2}(y)\int_y^\infty
p_{t/2}(x)\,dx\,dy\geq C_2\sqrt{t}\;e^{-\frac{\eps^2}{t}}.
\end{equation}
\end{lemma}

{\bf Proof:} By a change of variable and  (\ref{L1e}),
\begin{eqnarray}
&{}& \!\!\!\!\!\!\!\!\!\!\!\!\int_0^\infty p_{t/2}(y)\int_y^\infty
p_{t/2}(x)\,dx\,dy
\\
&\geq & \int_0^{\eps} p_{t/2}(y)\int_{\eps}^\infty
p_{t/2}(x)\,dx\,dy \nonumber \\
&=& \left (\int_0^{\eps}p_{t/2}(y)\,dy\right )\left (
\int_{\eps}^\infty p_{t/2}(x)\,dx\right )
\nonumber \\
&=& \left ( \frac{1}{\sqrt{\pi t}}\int_0^{\eps}e^{-y^2/t}\;dy\right
)\,\left ( \frac{1}{\sqrt{\pi
t}}\int_{\eps}^\infty e^{-x^2/t}\,dx\right )\nonumber \\
&=& C\,\left ( \int_0^{\frac{\eps}{\sqrt{t}}}e^{-\xi^2}\,d\xi\right
)\;\left (\int_{\frac{\eps}{\sqrt{t}}}^\infty
e^{-\xi^2}\,d\xi\right ) \nonumber \\
&\geq & C_2\sqrt{t}\;e^{-\frac{\eps^2}{t}}. \nonumber
\end{eqnarray}
\qed

\begin{lemma}\label{L3}
For any $\eps>0$, there exists a constant $C_3>0$ such that for every $t\in (0,1)$
and every $y\leq -(\aa+\eps)$,
\begin{equation}\label{L3e}
{\bf P}^y\left (X_s\leq -\aa,\;\;\forall s\in
[0,t/2]\right )\geq C_3.
\end{equation}
\end{lemma}

{\bf Proof:} By (\ref{p1}) we have
\begin{eqnarray}
&{}&\!\!\!\!\!\!\!\!\!\!\!\!{\bf P}^y\left (X_s\leq -\aa,\;\;\forall s\in
[0,t/2]\right )
\\
&\geq &{\bf P}^0(X_s\leq \eps,\;\forall s\in [0,1])\nonumber
\\
&=& 1-2{\bf P}^0(X_{1}\geq \eps)\nonumber = C_3>0
\end{eqnarray}
\qed

\begin{lemma}\label{L4}
For any $\eps>0$, there exist constants $\dd, C_4>0$ such that for every $t\in
(0,\dd)$,
\begin{equation}\label{L4e}
{\bf P}^0\left (X_s=0\;\;\;\hbox{for some}\;\;\; s\in
[t/2,t]\right )\geq C_4\;\sqrt{t}\;e^{-\frac{\eps^2}{t}}.
\end{equation}
\end{lemma}

{\bf Proof:} By the Markov property, the symmetry of Brownian
motion, (\ref{p1}), and Lemma \ref{L2},
\begin{eqnarray}
&{}&\!\!\!\!\!\!\!\!\!\!\!\!{\bf P}^0(X_s=0\;\;\;\hbox{for
some}\;\;\; s\in [t/2,t])
\\
&=&2\int_0^\infty p_{t/2}(y)\,{\bf P}^y(X_s=0\;\;\;\hbox{for
some}\;\;\; s\in [0,t/2])\,dy \nonumber
\\
&=&2\int_0^\infty p_{t/2}(y)\,{\bf P}^y(\tau_0\leq t/2)\,dy
\nonumber
\\
&=&2\int_0^\infty p_{t/2}(y)\,{\bf P}^0(\tau_y\leq t/2)\,dy
\nonumber
\\
&=& 4\int_0^\infty p_{t/2}(y)\,{\bf P}^0(X_{t/2}\geq
y)\,dy\nonumber
\\
&=&4\int_0^\infty p_{t/2}(y)\,\int_y^\infty
p_{t/2}(x)\,dxdy\nonumber\\
&\geq & C_4\sqrt{t}\;e^{-\frac{\eps^2}{t}}.\nonumber
\end{eqnarray}
\qed

\begin{lemma}\label{L5}
For any $\eps>0$, there exist constants $\dd, C_5>0$ such that for every $t\in
(0,\dd)$,
\begin{equation}\label{L5e}
{\bf P}^0\left (X_s\leq -\aa,\;\;\forall s\in [t/2,t]\right )\geq
C_5\;\sqrt{t}\;e^{-\frac{(\aa+\eps)^2}{t}}.
\end{equation}
\end{lemma}

{\bf Proof:} By the Markov property, Lemma \ref{L3}, and
(\ref{L1e}),
\begin{eqnarray}
&{}&\!\!\!\!\!\!\!\!\!\!\!\!{\bf P}^0(X_s\leq -\aa,\;\;\forall
s\in [t/2,t])
\\
&=&\int_{-\infty}^{-\aa}p_{t/2}(y)\;{\bf P}^y(X_s\leq
-\aa,\;\;\forall s\in [0,t/2])\;dy\nonumber
\\
&\geq & \int_{-\infty}^{-(\aa+\eps)}p_{t/2}(y)\;{\bf P}^y(X_s\leq
-\aa\;\;\forall s\in [0,t/2])\;dy\nonumber
\\
&\geq &
C\,\int_{-\infty}^{-(\aa+\eps)}p_{t/2}(y)\;dy=C\,\int_{(\aa+\eps)}^\infty
p_{t/2}(y)\;dy\nonumber
\\
&=&C\;\frac{1}{\sqrt{t}}\int_{(\aa+\eps)}^\infty
e^{-y^2/t}\;dy\nonumber
\\
&=& C\int_{\frac{(\aa+\eps)}{\sqrt{t}}}^\infty e^{-\xi^2}\;d\xi \nonumber
\\
&\geq &C \sqrt{t}\;e^{-\frac{(\aa+\eps)^2}{t}}.\nonumber
\end{eqnarray}
\qed

\bigskip

We can now prove Proposition \ref{prelim}. Fix $\aa$ with $0<\aa<
\frac{1}{\sqrt{2}}$ and choose $\eps>0$ so that $\eps^2 +
(\aa+\eps)^2 < \frac{1}{2} - \eps$. For this choice of $\eps$, fix
$\dd \in (0,1)$ appropriate for Lemmas \ref{L4}, and
\ref{L5}. Suppose $0<t<\dd$. We may apply Lemma \ref{L4} and Lemma \ref{L5}, using the independence of $X_s$ and $Y_s$, to get
\begin{eqnarray}\label{T2p3}
{\bf P}(T^{K_\aa}\leq t)&=& {\bf P}(Z_s\in (-\infty,-\aa],
\;\;\hbox{for some}\;\;\;s\in (0,t))  \\
&\geq & {\bf P}^0(X_s\leq -\aa, \;\;\hbox{for all}\;\;\;s\in
[t/2,t])  \nonumber
\\
&{}&\times \;\; {\bf P}^0(Y_s=0, \;\;\hbox{for some}\;\;\;s\in
[t/2,t])\nonumber
\\
&\geq & C\, t\,e^{-\frac{\eps^2}{t}}\;e^{-\frac{(\aa+\eps)^2}{t}}. \nonumber
\end{eqnarray}

On the other hand, it is proved in \cite{MAC} that for all $t>0$ and
all positive integers $n\geq 3$,
\begin{equation}\label{T2p4}
{\bf P}(T^\SD\leq t)\leq c(n)\;e^{-\frac{\cos^2(\pi/n)}{2t}}.
\end{equation}
Fixing $n$ large enough, we see that for all $t>0$,
\begin{equation}\label{T2p5}
{\bf P}(T^\SD\leq t)\leq C\;e^{-\frac{(\frac{1}{2}-\eps)}{t}}.
\end{equation}
Thus,
$$
\lim_{t \lar 0^+} \frac{{\bf P}(T^{K_\aa}<t)}{{\bf P}(T^\SD< t)}
\geq \lim_{t\to 0^+}\;C\;\frac{t\exp\left (
-\frac{\eps^2+(\aa+\eps)^2}{t}\right )}{\exp \left
(-(\frac{\frac{1}{2}-\eps}{t})\right )}=\lim_{t\to
0^+}\;t\;\exp\left (\frac{(\frac{1}{2}-\eps)-
(\eps^2+(\aa+\eps)^2)}{t}\right )=\infty.
$$
\qed

Now we are ready for the proof of Theorem \ref{T2}. By the scale
invariance of Brownian motion we can assume that $d(W) = 1$, and
rotation invariance allows us to assume that $-\aa \in \partial
U$, where $\aa=d(U) \in (0, \frac{1}{\sqrt{2}})$. Then clearly
$\SD \subseteq W$, and although it is not necessarily true that $U
\subseteq K_\aa$, we may still use our estimates for $K_\aa$ as a
lower bound by the following lemma.

\begin{lemma} \label{reflect}
$$
{\bf P}(T^U<t) \geq \frac{1}{2} {\bf P}(T^{K_\aa}< t).
$$
\end{lemma}

{\bf Proof:} Note that the complex conjugate of $Z_t$, $\bar Z_t$,
is also a Brownian motion. Let
$$
\tilde T_U = \inf\{t>0:\bar Z_t\notin U\}.
$$
We claim that $T_U \wedge \tilde T_U \leq T_{K_\aa}$ a.s. If not,
then the union of Brownian paths
$$
\{Z_t: 0 \leq t \leq T_{K_\aa}\}
\cup \{\bar Z_t: 0 \leq t \leq T_{K_\aa}\}
$$
would be a closed
curve separating $-\aa$ from $\ff$, and this contradicts simple
connectivity. Thus,
$$
{\bf P}(T_U \wedge \tilde T_U<t) \geq {\bf P}(T^{K_\aa}< t).
$$
But
$$
{\bf P}(T_U \wedge \tilde T_U<t) \leq {\bf P}(T_U <t) + {\bf
P}(\tilde T_U<t) = 2 {\bf P}(T_U <t),
$$
 and the lemma follows.
\qed

\bigskip

As for Theorem \ref{T2},
$$
\lim_{t \lar 0^+} \frac{{\bf P}(T^{U}<t)}{{\bf P}(T^W< t)} \geq
\lim_{t \lar 0^+} \frac{\frac{1}{2}{\bf P}(T^{K_\aa}<t)}{{\bf P}(T^\SD< t)}=\infty,
$$

completing the proof.

\medskip

{\bf Remark.} In fact, Lemma \ref{reflect} holds with the constant
$1$ in place of $\frac{1}{2}$. However, the proof requires a number of results on symmetrization and polarization which are not related to the rest of this paper; for this reason we have given the simpler result and proof above, and postpone the proof of the stronger result until the end of Section 3.

\subsection{Proof of Proposition \ref{T3}}

Let $p\in({\rm H}(W),{\rm H}(U))$ and $\delta=\frac{{\rm
H}(U)-p}{2}$. Since ${\bf E}[(T^{U})^{p+\frac{3}{2}\delta}]<\ff$,
the well-known Markov inequality (see e.g. \cite[6.17]{Fol})
implies
$$
{\bf P}(T^{U}>t)\leq\frac{{\bf E}[(T^{U})^{p+\frac{3}{2}\delta}]}{t^{p+\frac{3}{2}\delta}}.
$$
We now need a lower bound on ${\bf P}(T^{W}>t)$. For that purpose
we will use the so-called ``layer cake" representation for the
$p$-th moment (see e.g. \cite[6.24]{Fol}):
\begin{equation} \label{eq:1}
{\bf E}((T^{W})^{p})=p\int_{0}^{+\infty}t^{p-1}{\bf P}(T^{W}>t)dt.
\end{equation}

We now claim that
$$
\limsup_{t \to +\infty}\frac{{\bf P}(T^{W}>t)}{t^{-(p+\delta)}}
=+\infty,
$$
since otherwise $\frac{{\bf P}(T^{W}>t)}{t^{-(p+\delta)}}$ is
bounded above by a constant, and then by (\ref{eq:1}) we get ${\bf
E}[(T^{W})^{p}]<+\infty$ which contradicts the definition of ${\rm
H}(W)$. We obtain
\begin{eqnarray}
\limsup_{t \to +\infty} \frac{{\bf P}(T^{W}>t)}{{\bf P}(T^{U}>t)}
&\geq &  \limsup_{t \to +\infty}\frac{{\bf
P}(T^{W}>t)}{t^{-(p+\dd)}}\;\frac{t^{-(p+\dd)}\;t^{p+\frac{3}{2}\dd}}{{\bf
E}[(T^U)^{p+\frac{3}{2}\dd}}]\nonumber \\
&=& \limsup_{t \to +\infty} \frac{{\bf P}(T^{W}>t)}{t^{-(p+\dd)}}
\;\frac{t^{\frac{\dd}{2}}}{{\bf
E}[(T^{U})^{p+\frac{3}{2}\dd}]}=+\infty, \nonumber
\end{eqnarray}
which ends the proof. \qed

\section{Concluding remarks} \lll{conc}

As was discussed in the Introduction, Theorem $\ref{T2}$ shows
that the unit disk is not extremal among Schlicht domains for fast
exits. In fact, many Schlicht domains, including the Koebe domain
and the half-plane $\{{\rm Re}(z) \geq -\frac{1}{2}\}$, have a
higher probability of fast exits. We do not know if the constant
$\frac{1}{\sqrt{2}}$ in Theorem \ref{T2} is optimal; it would be
nice to know what is the best possible. As mentioned in the introduction, there is reason to suspect that the best possible constant is even 1. To see this, let $H_\aa = \{{\rm Re}(z) < \aa\}$. We then have the following proposition.

\begin{proposition} \label{}
 Let W be a simply connected domain
with $0\in W$. Suppose that $d(W) =1$, and $0 < \aa < 1$. Then, for all sufficiently small $t>0$,

\begin{equation}\label{T2e}
{\bf P}(T^{H_\aa}<t)> {\bf P}(T^W< t).
\end{equation}

In fact,

$$
\lim_{t \lar 0^+} \frac{{\bf P}(T^{H_\aa}<t)}{{\bf P}(T^W< t)} = \ff.
$$
\end{proposition}

{\bf Proof:} (sketch) By projection onto the real part, the exit time of $H_\aa$ has the same distribution as $\tau_\aa$, the first hitting time of $\aa$ by a one-dimensional Brownian motion, as defined at the beginning of Section \ref{2}. It is then straightforward to show using the Gaussian density that ${\bf P}(T^{H_\aa}<t) \geq C e^{\frac{-\aa^2(1+\eps)}{2t}}$ for $t$ sufficiently small and any $1> \eps >0$. On the other hand, $\SD \subseteq W$, and from \rrr{T2p5} we therefore have ${\bf P}(T^W\leq t)\leq C\;e^{-\frac{(1-\eps)}{2t}}$ for $t$ sufficiently small and any $1> \eps >0$. The result follows by choosing $\eps$ small enough so that $\aa^2 (1+\eps) < 1-\eps$. \qed

{\bf Remark:} Naturally, this result holds with $H_\aa$ replaced by any domain contained in a rotation of $H_\aa$.

\vski

The moments of the exit time have been considered previously by
several authors. In \cite{burk}, it is shown that for the wedge
$R_\th = \{|{\rm Arg}(z)|<\th\}$, that is, the infinite wedge
centered at the positive real axis of angular width $2\th$, we
have ${\bf E}[(T^{R_\th})^p] < \ff$ if and only if $p <
\frac{\pi}{4 \th}$, so
\begin{equation} \label{burkresult}
{\rm H}(R_\th) = \frac{\pi}{4\th}.
\end{equation}

A domain $W$ is {\it spiral-like of order $\si \geq 0$ with center $a$} if, for any $z \in W$, the spiral $\{a+(z-a) \mbox{ exp}(te^{-i \si}) : t \leq 0\}$ also lies within $W$; $W$ is {\it star-like} if it is spiral-like of order $\si = 0$. The quantity ${\rm H}(W)$ can be determined explicitly if $W$
is star-like or spiral-like, as is shown in \cite{markowsky}, with
equivalent analytic results appearing in \cite{hansen} and
\cite{hansenspi}. In particular, if we take $a=0$ then, since $W$ is spiral-like, the quantity

\begin{equation} \label{bigmax}
\AA_{r,W} = \max \{m(E): E \mbox{ is a subarc of } W \cap \{|z|=r\}\},
\end{equation}

is non-increasing in $r$ (here $m$ denotes angular Lebesgue measure on the circle). We may therefore let $\AA_W = \lim_{r \nearrow \ff} \AA_{r,W}$, and then \cite[Thm. 2]{markowsky} we have $H(W) = \frac{\pi}{2 \AA_W \cos^2 \si}$.

\vski

As mentioned before we suspect that in many cases the $\limsup$ in Proposition \ref{T3} is not necessary, and venture the following conjecture.

\begin{conject} \label{Con1}
Suppose that ${\rm H}(U) > {\rm H}(W)$ and $W$ is spiral-like.
Then
\begin{equation}\label{T3ee}
\lim_{t \lar \ff} \frac{{\bf P}(T^W > t)}{{\bf P}(T^U > t)} = \ff.
\end{equation}
\end{conject}

Note that this includes the case that $W$ is star-like, as well as the wedge $R_\th$.
This conjecture would follow from the following, if true.

\begin{conject} \label{}
Suppose that $W$ is spiral-like. Then for any $p>{\rm H}(W)$ there
is a constant $C>0$ so that ${\bf P}(T_W>t) \geq \frac{C}{t^p}$.
\end{conject}

Our evidence for the truth of Conjecture \ref{Con1} is at follows. First
note that Markov's inequality, used as in Proposition
\ref{T3}, yields the following fact.

\begin{proposition} \label{-}
For any $p<{\rm H}(U)$, there is a constant $C>0$ so that ${\bf
P}(T^U>t) \leq \frac{C}{t^p}$.
\end{proposition}

Furthermore the bound required is true in the $\limsup$ sense, as
is shown in the proof of Proposition \ref{T3}. Next we prove
Conjecture \ref{Con1} when $W$ is a half-plane or quarter-plane. Let $W =
\{{\rm Re}(z) < 1\}$; recall from \rrr{burkresult} that ${\rm
H}(W) = \frac{1}{2}$. Then, using the reflection principle, ${\bf
P}(T^U>t)$ can be bounded below as follows.

\begin{eqnarray}\label{L8p1}
&{}&\!\!\!\!\!\!\!\!\!\!\!\!{\bf P}^0(X_s<1,\;\;\;\forall s\in
[0,t])\\
&=&1-2{\bf P}^0(X_t\geq 1)=1-2\frac{1}{\sqrt{2\pi
t}}\int_{1}^\infty e^{-x^2/2t}\;dx\nonumber \\
&=&\frac{2}{\sqrt{\pi}}\int_0^{\frac{1}{\sqrt{2t}}}
e^{-\xi^2}\;d\xi\geq \frac{C}{\sqrt{t}}.\nonumber
\end{eqnarray}

Now suppose $W = \{{\rm Re}(z) > 0, {\rm Im}(z)>0\}$; recall from
\rrr{burkresult} that ${\rm H}(W) = 1$. Then, using the
independence of the one dimensional components of planar Brownian
motion, and the calculation for the half-plane, we have

\begin{eqnarray}\label{T3p3}
{\bf P}^{1+i}(T^W>t) &= & {\bf P}^1(X_s>0,\;\;\;\forall s\in
[0,t])^2
\geq \frac{C}{t}\nonumber
\end{eqnarray}

for $t$ bounded away from $0$. Finally, we prove an improvement of Lemma \ref{reflect}, as we
promised in Section 2.

\begin{proposition} \label{symmetrization}
Let $U$ be a simply connected domain with $0\in U$ and let
$\aa=d(U)\in (0,\infty)$. If $K_\aa=\SC\setminus (-\infty,-\aa]$,
then for every $t>0$,
\begin{equation}\label{syme}
{\bf P}(T^U>t)\leq {\bf P}(T^{K_\aa}>t).
\end{equation}
\end{proposition}

{\bf Proof:} Let $p^U(t,0,w)$ be the transition density function
for Brownian motion killed upon hitting $\partial U$. Then (see
e.g. \cite[Theorem 2.4]{CZ})
\begin{equation}\label{sym1}
{\bf P}(T^U>t)=\int_U p^U(t,0,w)\;A(dw),\;\;\;0<t<+\infty,
\end{equation}
where $A$ denotes the area measure. A similar formula holds for
${\bf P}(T^{K_\aa}>t)$. Therefore, it suffices to prove that
\begin{equation}\label{sym2}
\int_U p^U(t,0,w)\;A(dw)\leq \int_U p^{K_\aa}(t,0,w)\;A(dw),
\;\;\;0<t<+\infty.
\end{equation}

\medskip

We will prove (\ref{sym2}) using the theory of polarization and
symmetrization. We refer to \cite{Hay}, \cite{Dub}, \cite{BS} for
the definitions and basic facts.

\vski

The function $p^U(t,0,w)$ satisfies the heat equation on $U$. Let
$P_HU$ denote the polarization of $U$ with respect to a half-plane
$H$ with $0\in\partial H$. Then \cite[Theorem 9.4]{BS}
\begin{equation}\label{sym3}
p^U(t,0,w)+p^U(t,0,R_Hw)\leq
p^{P_HU}(t,0,w)+p^{P_HU}(t,0,R_Hw),\;\;\;0<t<\infty,
\end{equation}
where $R_Hw$ denotes the reflection of $w$ in the line $\partial
H$. It follows from (\ref{sym3}) that for every $r\in
(0,+\infty)$,
\begin{equation}\label{sym4}
\int_0^{2\pi}p^U(t,0,re^{i\theta})d\theta\leq
\int_0^{2\pi}p^{P_HU}(t,0,re^{i\theta})d\theta,\;\;\;0<t<\infty.
\end{equation}

By a standard technique involving a sequence of polarizations,
(\ref{sym4}) leads to the inequality
\begin{equation}\label{sym5}
\int_0^{2\pi}p^U(t,0,re^{i\theta})d\theta\leq
\int_0^{2\pi}p^{U^*}(t,0,re^{i\theta})d\theta,\;\;\;0<t<\infty,\;0<r<\infty,
\end{equation}
where $U^*$ is the circular symmetrization of $U$ with respect to
the positive semi-axis. Since $U^*\subset K_\aa$, we have
\begin{equation}\label{sym6}
\int_0^{2\pi}p^{U^*}(t,0,re^{i\theta})d\theta\leq
\int_0^{2\pi}p^{K_\aa}(t,0,re^{i\theta})d\theta,\;\;\;0<t<\infty,\;0<r<\infty.
\end{equation}
By (\ref{sym5}), (\ref{sym6}), and integration over $r\in
(0,\infty)$, we obtain (\ref{sym2}).
 \qed



\bibliographystyle{alpha}
\bibliography{CAbib}

\end{document}